\documentclass{amsart}
\usepackage{amsmath}
\usepackage{amssymb}
\newtheorem{Lemma}{Lemma}[section]
\newtheorem{Th}[Lemma]{Theorem}
\newtheorem{Prop}[Lemma]{Proposition}

\newtheorem{Cor}[Lemma]{Corollary}

\newtheorem{Remark}[Lemma]{Remark}
\newenvironment{Proof}{{\sc Proof.}\ }{~\rule{1ex}{1ex}\vspace{0.2truecm}}
\newcommand{\Cal}[1]{{\mathcal #1}}
\newcommand{\End}{\mbox{\rm End}}

\begin{document}
    \title[Rings with right perfect factors]{Rings whose proper factors are right perfect}
    \author[Alberto Facchini and Catia Parolin]{Alberto Facchini and Catia Parolin}
    \address{Dipartimento di Matematica Pura e Applicata, Universit\`a di Padova, 35121 Padova, Italy}
    \email{facchini@math.unipd.it, catia.parolin@studenti.unipd.it}
\thanks{Partially
supported by Ministero dell'Istruzione, dell'Universit\`a e della Ricerca, Italy (Prin
2007 ``Rings, algebras, modules and categories'')  and by  Universit\`a di Padova (Progetto di Ricerca di Ateneo CPDA071244/07).}

 \begin{abstract} We show that practically all the properties of almost perfect rings discovered by Bazzoni and Salce in the paper \cite{salceapd}  for commutative rings hold in the non-commutative setting. 
 \end{abstract}
    \maketitle

\section{Introduction}

We say that a ring $R$ is {\em right almost perfect} if $R/I$ is a right perfect ring for every proper non-zero two-sided ideal $I$ of $R$. Similarly we define \emph{left almost perfect} rings. Commutative almost perfect rings were defined by Bazzoni and Salce in \cite{JLMS}. The class of commutative almost perfect rings has several interesting properties and characterization, and has been studied by Bazzoni, Fuchs, Sang Bum Lee, Salce, Zanardo and others. We find that most of the properties of these rings presented in \cite{salceapd} still hold in the non-commutative setting: non-prime right almost perfect rings are right perfect, over the prime rings the non-faithful modules are semiartianian, etcetera. Our main results are Theorems~\ref{finalmente} and  \ref{hloc}. 

We also introduce the notion of (non-commutative) {\em $h$-local} ring. This also is an extension to non-commutative rings of the corresponding notion studied  in the commutative setting first by Jaffard \cite[Th.~6]{Jaffard}, and then by Matlis \cite{Cotorsion} and others. In Proposition \ref{noeth}, we show that one of the characterizations of commutative almost perfect ring holds for the right noetherian right almost perfect rings, but not for the non-noetherian ones. 

Bazzoni and Salce also proved that a commutative integral domain $R$
is almost perfect if and only if it  is $h$-local and every localization of $R$ at a maximal ideal is almost perfect.
As localization at maximal ideals is typical of the commutative setting, this property clearly does not have an obvious non-commutative counterpart. The next step now is to see whether the other properties of commutative almost perfect rings, for instance those in \cite{JLMS}, also hold in the non-commutative case.

Some examples are given.
In this paper, all rings $R$ are associative rings with identity $1\ne 0$, and $J(R)$ denotes the Jacobson radical of the ring $R$.

\section{Bazzoni and Salce's results}

We will very frequently use in the sequel the various characterizations of right perfect rings that appear in the celebrated result due to Bass and called {\em Theorem~P} \cite[Theorem \ 28.4]{andersonfuller}. Hence we state it here for later reference. Recall that a (non-necessarily commutative) ring $R$ is said to be: (1) {\em semilocal} if $R/J(R)$ is semisimple artinian; (2) a \emph{right max ring} if every non-zero right $R$-module has a maximal submodule; and (3) {\em right semiartinian} if every non-zero right $R$-module has a minimal submodule (equivalently, if every right $R$-module is an essential extension of its socle).

\begin{Th}\label{28.4} {\rm (Bass's Theorem P)}
The following conditions are equivalent for a ring $R$:
\begin{enumerate}
  \item $R$ is right perfect;
  \item $R$ is semilocal and its Jacobson radical $J(R)$ is right $T$-nilpotent;
  \item $R$ is semilocal and right max;
  \item every flat right $R$-module is projective;
  \item $R$ satisfies the descending chain condition on principal left ideals;
  \item $R$ is left semiartinian and contains no infinite orthogonal set of idempotents.
\end{enumerate}
\end{Th}

We now briefly recall the main results on commutative almost perfect rings obtained by Bazzoni and Salce in \cite{salceapd}. They defined {\em almost perfect rings} as those commutative rings $R$ for which $R/I$ is a perfect ring for every proper non-zero ideal $I$ of $R$. Recall that a commutative domain $R$ is \emph{$h$-local} \cite{Cotorsion} if  $R/I$ is semilocal for every non-zero ideal $I$ of $R$ and $R/P$ is local for every non-zero prime ideal $P$ of $R$.
In the commutative case, a ring $R$ is perfect if and only if it is semilocal and every localization of $R$ at a maximal ideal is a perfect ring \cite[Theorem \ 1.1]{salceapd}. The following are the four main results contained in the paper by Bazzoni and Salce \cite{salceapd}.

\begin{Prop}\label{s1.3} {\rm \cite[Proposition~1.3]{salceapd}}
If $R$ is an almost perfect commutative ring and $R$ is not an integral domain, then $R$ is perfect.
\end{Prop}

\begin{Th}\label{s2.2} {\rm \cite[Theorem~4.4.1]{enochsjenda}}
The following conditions are equivalent  for a commutative integral domain $R$: \begin{enumerate}
 \item every non-zero torsion module contains a simple module;
 \item every torsion $R$-module is semiartinian;
 \item for every non-zero ideal $I$ of $R$, $R/I$ contains a simple module;
 \item for every non-zero $R$-submodule $A$ of $Q$, $Q/A$ contains a simple module;
 \item $Q/R$ is semiartinian.
\end{enumerate}
\end{Th}

\begin{Th}\label{s2.3} {\rm \cite[Theorem~2.3]{salceapd}}
 The following conditions are equivalent for a commutative integral domain $R$:
\begin{enumerate}
 \item $R$ is almost perfect;
 \item $R$ is $h$-local and every localization of $R$ at a maximal ideal is almost perfect; 
 \item $R$ is $h$-local and satisfies one of the equivalent conditions of Theorem~\ref{s2.2}.
\end{enumerate}
\end{Th}

\begin{Cor}
 The following conditions are equivalent for a commutative local integral domain $R$:\begin{enumerate}
\item $R$ is almost perfect;
\item $Q/R$ is semiartinian;
\item every non-zero torsion module is semiartinian.
\end{enumerate}
\end{Cor}

\section{Non-commutative almost perfect rings}

We will now see how the previous results modify when the ring $R$ is not commutative. As we have already said in the introduction, we call a ring $R$ {\em right almost perfect} if $R/I$ is a right perfect ring for every proper ideal $I\ne 0$ of $R$.

We begin with some immediate examples of right almost perfect rings. 

(1) Right perfect rings are right almost perfect \cite[Corollary~28.7]{andersonfuller}. 

(2) Simple rings are trivially right and left almost perfect. Another very interesting example is given by the nearly simple chain rings due to Dubrovin and considered by Puninski \cite{dubrovin, puninski1, puninski2}. 
A {\em nearly simple} chain ring $R$ is a non-commutative right and left chain ring (i.e., the right ideals and the left ideals are linearly ordered under inclusion), with exactly three two-sided ideals, necessarily the ideals $0$, $R$ and the maximal ideal $J(R)$. As chain rings $R$ are local, $R/J(R)$ is necessarily a division ring, hence a perfect ring. This shows that nearly simple chain rings are right and left almost perfect. Puninski considers both examples of nearly simple chain rings that are (non-commutative) integral domains, and 
examples of nearly simple chain rings that are prime rings but not integral domains. 

(3) There are left almost perfect rings that are not right almost perfect. For instance, let $k$ be a field and let $k_{\omega}$ be the $k$-algebra of all matrices with entries in $k$, countably many rows and columns, and in which each row has only finitely many non-zero entries. Let $N$ be the set of all strictly lower triangular matrices in $k_{\omega}$ with only finitely many non-zero entries and let $R$ be the subalgebra $k+N$ of $k_{\omega}$. If we denote by $E_{i,j}$ the unitary matrices, where $i$ and $j$ are positive integers, one has that $E_{i,j}\in N$ if and only if $i>j$. It is known that $R$ is left perfect but not right perfect \cite[Example (5), p.~476]{Bass}. The Jacobson radical of $R$ is $N$. In particular, $R$ is left almost perfect. In order to show that $R$ is not right almost perfect, consider the principal two-sided ideal $I$ of $R$ generated by $E_{2,1}$. It is easily seen that $I$ is the vector space generated by all products $E_{i,j}E_{2,1}E_{k,\ell}$, with $i>j$ and $k>\ell$, that is, $I$ is the vector space generated by all $E_{i,1}$ with $i\ge 2$. Then $I\leq J(R)$, so that $J(R/I)=J(R)/I=N/I$. In the sequence $\dots, E_{5,4},E_{4,3},E_{3,2}$, all the elements are in $N$ and the products  $E_{n,n-1}E_{n-1,n-2}\dots E_{3,2}=E_{n,2}$ are not in $I$. This proves that $N/I$ is not right $T$-nilpotent, and so $R/I$ is not right perfect. This shows that $R$ is left almost perfect but  not right almost perfect. 

(4) There is no relation between being an almost perfect ring and being a semiperfect ring. For instance, the ring $\mathbb{Z}$ of integers is almost perfect but not semiperfect, and a commutative valuation domain of Krull dimension two is semiperfect but not almost perfect. 

Our first result is the non-commutative analogue of Proposition~\ref{s1.3}. It shows that, for non-prime rings, the notion of right almost perfect ring and right perfect ring coincide.

\begin{Th}\label{finalmente} If a ring $R$ is right almost perfect and not a prime ring, then $R$ is right perfect.\end{Th}

\begin{Proof} Let $R$ be a right almost perfect ring that is not a prime ring. We will distinguish two cases.

{\em First case:} $R$ has a non-zero nilpotent two-sided ideal.

In this case, $R$ has a two-sided ideal $K\ne 0$ with $K^2=0$. In particular, $K$ is nilpotent, hence $K\subseteq J(R)$. Thus $J(R)\ne 0$, so that $R/J(R)$ is right perfect, hence semisimple artinian (Theorem~\ref{28.4}). Thus $R$ is semilocal, and so it has no infinite orthogonal set of idempotents. In order to conclude, by Theorem~\ref{28.4}, it suffices to show that every non-zero left $R$-module contains a simple submodule. Let $_RM\ne 0$ be a  left $R$-module. If $KM=0$, then $M$ is a left $R/K$-module. But $R/K$ is right perfect, so that $M$ has a simple $R/K$-submodule, which is also a simple $R$-submodule. If $KM\ne 0$, then $KM$ is a non-zero left  $R/K$-module. But $R/K$ is right perfect, so that $KM$ has a simple $R/K$-submodule, which is also a simple $R$-submodule. Thus $M$ has a simple $R$-submodule. Thus $R$ is right perfect in this first case.

{\em Second case:} $R$ has no non-zero nilpotent two-sided ideals.

Since $R$ is not a prime ring, $R$ has two non-zero two-sided ideals $I$ and $J$ with $IJ=0$. Then $(I\cap J)^2\subseteq IJ=0$, so that $I\cap J=0$ because $R$ has no non-zero nilpotent two-sided ideals. We will now show that $R$ contains no infinite orthogonal set of idempotents. Assume the contrary, and let $E$ be an infinite orthogonal set of distinct idempotents of $R$. Then $E_I:=\{\, e+I\mid e\in E\,\}$ is an orthogonal set of idempotents of $R/I$, which is right perfect, so that $E_I$ must be a finite set. It follows that there is a partition of $E$ into finitely many subset $E_1,\dots,E_n$ with the property that, for every $e,f\in E$, $e-f\in I$ if and only if $e$ and $f$ belong to the same block $E_i$ of the partition. Since $E$ is infinite, one of the blocks, $E_t$ say, is an infinite set. Thus $E_t$ is an infinite orthogonal set of distinct idempotents of $R$. The set $E_{t,J}:=\{\, e+J\mid e\in E_t\,\}$ is an orthogonal set of idempotents of $R/J$, which is right perfect, so that $E_{t,J}$ must be a finite set. It follows that there is a partition of $E_t$ into finitely many subset $E'_1,\dots,E_m'$ with the property that, for every $e,f\in E_t$, $e-f\in J$ if and only if $e$ and $f$ belong to the same block $E'_j$ of this partition of $E_t$. As $E_t$ is infinite, one of these blocks, $E'_\ell$ say, is infinite. But for every $e,f\in E'_\ell$, $e-f\in J$ because $e$ and $f$ belong to the same block $E'_\ell$, and $e-f\in I$ because both $e$ and $f$ belong to $E_t$. Thus $e-f\in I\cap J=0$, i.e., $e=f$ for every $e,f\in E'_\ell$. This shows that $E'_\ell$ has exactly one element, which contradicts what we had previously proved. This proves that $R$ contains no infinite orthogonal set of idempotents. In order to conclude, by Theorem~\ref{28.4}, it suffices to show that every non-zero left $R$-module contains a simple submodule. Let $_RM\ne 0$ be a  left $R$-module. If $IM=0$, then $M$ is a left $R/I$-module. But $R/I$ is right perfect, so that $M$ has a simple $R/I$-submodule, hence a simple $R$-submodule. If $IM\ne 0$, then $IM$ is a non-zero left  $R/J$-module, because $(JI)^2=J(IJ)J=0$ implies $JI=0$ ($R$ has no non-zero nilpotent two-sided ideals). But $R/J$ is right perfect, so that $IM$ has a simple $R/J$-submodule, hence a simple $R$-submodule. Thus $M$ has a simple $R$-submodule. This proves that $R$ is right perfect in this second case also.\end{Proof}

Puninski's example of a nearly simple chain domain $R$ \cite{puninski1} immediately shows that the conditions of Theorem~\ref{s2.2} do not hold for the right almost perfect ring~$R$. That ring $R$ is an Ore domain, hence it has its natural torsion theory. In order to recall the definition of this ring $R$, let $G$ be the group of affine linear functions on~$\mathbb{Q}$:
$$G=\{\,\alpha_{a,b}\colon \mathbb{Q}\to \mathbb{Q}\mid  a,b\in \mathbb{Q},\: a>0\,\},
$$
where $\alpha_{a,b}(t)=at+b$ and the group operation on $G$ is the composition of functions. Fix a a positive irrational number $\epsilon$ in the field of real numbers.
It is possible to define a right order on $G$ with generalized positive cone $P:=\{\,\alpha_{a,b}\in G\mid \epsilon\leq \alpha_{a,b}(\epsilon)\,\}$. Let $k$ be a division ring, $k[P]$ the semigroup ring, and consider $M:=\sum_{\alpha\in P^+}\alpha k[P]$ in $k[P]$, where $P^+=\{\,\alpha_{a,b}\in P\mid \epsilon <\alpha_{a,b}(\epsilon)\,\}$. Then $M$ is a maximal ideal in $k[P]$. The subset $k[P]\setminus M$ is a right and left Ore set in $k[P]$. Let $R$ denote the localization of $k[P]$ with respect to $k[P]\setminus M$. The ring $R$ is a nearly simple chain domain, hence it is a right and left almost perfect ring. Moreover, $R$ is a local ring and an Ore domain. Thus the natural torsion theory for this ring is the usual torsion theory. 
The valuation group of $R$ is $G$. This group is not a discrete group. It follows that $Q/R$, where $Q$ is the division ring of fractions of $R$, is a torsion $R$-module with zero socle, so that the conditions of Theorem~\ref{s2.2} do not hold for the almost perfect ring ring $R$.

\medskip

Let us pass to consider the structure of right perfect rings. We begin with an elementary lemma.

\begin{Lemma}\label{prime}
Every prime right perfect ring is a simple artinian ring.
\end{Lemma}

\begin{Proof}
If $R$ is a prime right perfect ring, then $R$ satisfies the descending chain condition on principal left ideals by Theorem~\ref{28.4}, so that $R$ is semisimple artinian by \cite[Theorem~10.24]{Lam}. Since $R$ is prime, it must be a simple ring.\end{Proof}

\begin{Cor} A non-zero two-sided ideal of a right almost perfect ring $R$ is a maximal ideal if and only if it is a prime ideal, if and only if it is a right primitive ideal.\end{Cor}

\begin{Proof} Every maximal ideal is prime. If $I\ne 0$ is a prime ideal, then $R/I$ is a prime right perfect ring, so that $R/I$ is simple artinian by the previous Lemma. Thus $I$ is the right annihilator of the unique simple right $R/I$-module. Finally, let $I\ne 0$ be a right primitive ideal. Then $R/I$ is a right perfect ring with a faithful simple right module. In particular, $J(R/I)=0$. But right perfect rings are semisimple artinian modulo their Jacobson radical, so that $R/I$ is a semisimple artinian ring. Since it has a faithful simple right module, it follows that $R/I$ is simple artinian. Thus $I$ is maximal in $R$.\end{Proof}

We say that a prime ring $R$ is \emph{$h$-local} if: (1) for every non-zero proper two-sided ideal $I$ of $R$, the factor ring $R/I$ is semilocal; and (2) every non-zero prime two-sided ideal of $R$ is contained in only one maximal two-sided ideal of $R$.

Clearly, local prime rings are $h$-local. From Theorem~\ref{28.4} and Lemma~\ref{prime}, we immediately get that:

\begin{Cor}\label{aphsemiloc}
 If $R$ a prime right almost perfect ring, then $R$ is $h$-local.
\end{Cor}

Any prime ring $R$ is a right linearly topological ring in a natural way \cite[p.~144]{stenstrom}. Let $\Cal B$ be the set of all non-zero two-sided ideals of $R$. The topology on the prime ring $R$  has $\Cal B$ as a basis of neighborhoods of $0$. It is not a Hausdorff topology in general. For instance, if $R$ is a nearly simple prime chain ring, the closure of zero in this topology is the maximal ideal of $R$. Moreover, this right ideals of the prime ring $R$ that are open in this right linear topology do not form a right Gabriel topology in general, but only a {\em divisible right Oka family} \cite{Reyes}. More precisely, recall that if $I$ is a right ideal of a ring $R$, the {\em core} of $I$ is the largest two-sided ideal of $R$ contained in $I$. It coincides with the annihilator in $R$ of the right $R$-module $R/I$. If $R$ is a prime ring, the class $\Cal T$ of all right $R$-modules whose elements are annihilated by an element of $\Cal B$ is closed under submodules, homomorphic images and direct sums. 

Non-faithful modules are sometimes called {\em bounded}. From Theorem~\ref{28.4} (Condition (6)) and the standard characterizations of right semiartinian rings \cite[Proposition~VIII.2.5]{stenstrom}, we immediately get the next two results. Proposition~\ref{equivalenze} gives the non-commutative analogue of Theorem~\ref{s2.2}.

\begin{Lemma}\label{apsemiart}
 If $R$ is a right almost perfect ring, then every non-faithful left $R$-module is semiartinian.
\end{Lemma}

\begin{Prop}\label{equivalenze}
The following conditions are equivalent for a ring $R$:
\begin{enumerate}
  \item{every non-faithful left $R$-module is semiartinian;}
  \item{every non-faithful left $R$-module is an essential extension of its socle;}
  \item{every non-zero non-faithful left $R$-module contains a simple submodule;}
  \item{$R/I$ is left semiartinian for every non-zero two-sided ideal $I$ of $R$;}
  \item{for every proper left ideal $L$ with a non-zero core, the cyclic left $R$-module $R/L$ contains a simple submodule.}
\end{enumerate}
\end{Prop}

%


%

We are ready to state and prove the non-commutative analogue of Theorem~\ref{s2.3}

\begin{Th}\label{hloc}
The following statements are equivalent for a ring $R$:
\begin{enumerate}
  \item{the ring $R$ is right almost perfect;}
  \item{the ring $R$ is $h$-local and it satisfies one of the equivalent condition of Proposition~\ref{equivalenze}.}
 \end{enumerate}
\end{Th}

\begin{Proof}
  (1)$\Leftrightarrow$(2). We have already shown that (1) implies that $R$ is $h$-local and every non-faithful left $R$-module is semiartinian (Corollary~\ref{aphsemiloc} and Lemma~\ref{apsemiart}).\\
  Let $R$ be an $h$-local ring and suppose that every non-faithful left $R$-module is semiartinian. We want to show that $R$ is right almost perfect. Let $I$ be a proper non-zero two-sided ideal of $R$. Trivially every non-zero left $R/I$-module has a simple submodule (because every non-faithful left $R$-module is semiartinian and so it contains a simple submodule).
 The ring $R/I$ is semilocal, hence has no infinte orthogonal set of idempotents. By \cite[Theorem~{28.4}(f)]{andersonfuller}, this implies that $R/I$ is right perfect.
\end{Proof}

\begin{Cor}
Let $R$ be a local ring. Then $R$ right almost perfect if and only if it satisfies one of the equivalent conditions of Proposition~\ref{equivalenze}.
\end{Cor}

\section{Right noetherian rings}

For right noetherian prime rings, the situation is much better than in the previous section, in which the prime ring $R$ was arbitrary. Recall that a right noetherian ring is right perfect if and only if it is right artinian (To see this, let $R$ be a right noetherian right perfect ring. Then $R$ right perfect implies $R$ semilocal and $J(R)$ right $T$-nilpotent. In particular, $J(R)$ is nil, hence nilpotent by Levitzki's Theorem \cite[Theorem~15.22]{andersonfuller}. By Hopkins's Theorem \cite[Theorem~15.20]{andersonfuller}, $R$ is right artinian.)

Let $R$ be a right noetherian prime ring. Let $\Cal B$ be the set of all non-zero two-sided ideals of $R$. By \cite[Proposition~6.10]{stenstrom}, the set $\Cal B$ now {\em is} a basis for a right Gabriel topology. We will call {\em torsion $R$-modules} the right $R$-modules that are in the corresponding torsion class $\Cal T$. 

\begin{Prop}\label{noeth} The following conditions are equivalent for a right noetherian prime ring $R$:

(a) Every non-zero torsion right $R$-module contains a simple module.

(b) Every torsion right $R$-module is semiartinian.

(c) For every non-zero proper two-sided ideal $I$ of $R$, the ring $R/I$ is right semiartinian.

(d) For every non-zero proper two-sided ideal $I$ of $R$, the ring $R/I$ is right artinian.

(e) $R$ is a right almost perfect ring.
\end{Prop}

\begin{Proof} (a)${}\Rightarrow{}$(b) follows from the fact that every homomorphic image of a torsion module is torsion and the definition of semiartinian module.

(b)${}\Rightarrow{}$(c). As $I$ is non-zero, every right $R/I$-module is a torsion $R$-module, hence is a semiartinian $R$-module. It follows that it is semiartinian as an $R/I$-module also.

(c)${}\Rightarrow{}$(d). Since $R$ is right noetherian, the Gabriel topology corresponding to semiartinian modules consists of the right ideal $A$ of $R$ with $R/A$ of finite length. If $I$ is two-sided and $R/I$ is right semiartinian, then $1+I$ is annihilated by $I$, hence $I$ belongs to that Gabriel topology, and so $R/I$ is a right $R$-module of finite length. Thus $R/I$ is a right artinian ring.

(d)${}\Rightarrow{}$(e)  follows immediately from the fact that every right artinian ring is right perfect.

(e)${}\Rightarrow{}$(a).  If $M_R$ is a non-zero torsion module over the right almost perfect ring $R$, then $M_R$ has a non-zero element $x$, and $Rx$ is non-faithful because it is annihilated by the two-sided ideal that annihilates $x$. Therefore $Rx$, hence $M_R$, contains a simple module.
\end{Proof}


This generalizes Bazzoni-Salce in the case of right noetherian non-necessarily commutative rings.

\begin{Remark} {\rm One could think that over an h-local ring $R$, every torsion module $M_R$ decomposes into its primary components, one for each maximal ideal, like for commutative h-local rings, but this is not true, not even in the best of the cases we are dealing with. In our example the ring $R$ will be right and left artinian, hence right and left noetherian, right and left perfect, and semilocal. In particular, it is h-local. Let $k$ be a field and $R:=\binom{k\ k}{0\ k}$ the ring of all upper triangular $2\times 2$ matrices, so that $R$ has all the properties mentioned in the previous sentence. This ring $R$ has exactly two maximal two-sided ideals. The cyclic projective module $M_R:=\binom{1\ 0}{0\ 0}R=\binom{k\ k}{0\ 0}$ has its endomorphism ring isomorphic to $k$, hence is an indecomposable $R$-module, and is the extension of two non-isomorphic simple modules. Therefore $M_R$ does not decompose into two primary components corresponding to the two simple right $R$-modules.}\end{Remark}

We conclude the paper with a further remark. There are two mutually exclusive cases according to if the ideal 0 is maximal or not maximal. Also, a ring can be semilocal or not. Correspondingly, right almost perfect rings belong to exactly one of the following three classes:

First class: simple rings. It corresponds to the case in which $0$ is a maximal ideal. We already know that simple rings are right almost perfect.

Second class: it corresponds to the case in which $0$ is not a maximal ideal and $R$ is semilocal. A ring $R$ belongs to this class if and only if it is a semilocal non-simple ring, the factor rings $R/P$ are simple artinian for all non-zero prime ideals $P$ of $R$ and every non-zero non-faithful left $R$-module is semiartinian. Notice that these conditions are either two-sided conditions on the ring or conditions on {\em left} $R$-modules, and they characterize semilocal non-simple {\rm right} almost perfect rings. For instance, nearly simple chain rings belong to this class.

Third class: it corresponds to the case in which $0$ is not a maximal ideal and $R$ is not semilocal. A ring $R$ belongs to this class if and only if it is a non-simple ring, $J(R)=0$, every non-zero element of $R$ belongs to only finitely many maximal ideals of $R$, the factor rings $R/P$ are simple artinian for all non-zero prime ideals $P$ of $R$ and every non-zero non-faithful left $R$-module is semiartinian. 

\bigskip

We conclude with an example of a ring which belongs to this third class. It is due to Faith and Michler-Villamayor \cite[Remark~4.5]{MV}. It is a Von Neumann regular ring that is a right $V$-ring that is not a left $V$-ring. Recall that right $V$-rings are right max rings \cite[Theorem~2.1]{MV}.

\smallskip

Let $k$ be a field and let $V_k$ be an infinite dimensional vector space over $k$. Consider the endomorphism ring $E:=\End(V_k)$ and its two-sided ideal $S$ consisting of endomorphisms of finite rank. It is easily seen that the $k$-subalgebra $R:=k+S$ of $E$ has just three two-sided ideals: $0$, $R$ and $S$. In particular,  $R$ is prime. Trivially, the ring $R$ is right and left almost perfect. The Jacobson radical $J(R)$ of $R$ is zero. In fact, $J(R)$ is a proper two-sided ideal, so it can only be $0$ or $S$; but the element $1-E_{11}\in 1+S$ is not invertible. Thus $J(R)=0$.  In particular $R$ is not semilocal, otherwise $R/J(R)\cong R$ would be semisimple artinian, which is not. Finally $R$ not semilocal implies that $R$ is not right perfect and is not left perfect.

\end{document}